\documentstyle[11pt]{article}
\def\bkR{{\rm I\kern-.17em R}}
\def\bkC{{\rm \kern.24em \vrule width.05em height1.4ex depth-.05ex \kern-.26em C}}
\def\bkN{{\rm \kern.50em \vrule width.05em height1.4ex depth-.05ex \kern-.26em N}}
\def\D{{\cal D}}

\begin{document}

\author{Nuno Costa Dias\footnote{{\it ncdias@meo.pt}} \\ Jo\~{a}o Nuno Prata\footnote{{\it joao.prata@mail.telepac.pt}} \\ {\it Departamento de Matem\'atica} \\
{\it Universidade Lus\'ofona de Humanidades e Tecnologias} \\ {\it Av. Campo Grande, 376, 1749-024 Lisboa, Portugal}\\
{\it and}\\
{\it Grupo de F\'{\i}sica Matem\'atica}\\
{\it Universidade de Lisboa}\\
{\it Av. Prof. Gama Pinto 2}\\
{\it 1649-003 Lisboa, Potugal}}
\date{}
\title{A multiplicative product of distributions and a global formulation of the confined Schr\"odinger equation\footnote{Presented by N.C. Dias at the $3^{rd}$ Baltic-Nordic Workshop "Algebra, Geometry, and Mathematical Physics", G\"oteborg, Sweden, October 11-13, 2007.}}

\maketitle

\begin{abstract}
A concise derivation of a new multiplicative product of Schwartz distributions is presented.
The new product $\star$ is defined in the vector space ${\cal A}$ of piecewise smooth functions on $\bkR$ and all their (distributional) derivatives; it is associative, satisfies the Leibnitz rule and reproduces the usual product of functions for regular distributions. The algebra $({\cal A},+,\star)$ yields a sufficiently general setting to address some interesting problems. As an application we consider the problem of deriving a global formulation for quantum confined systems.
\end{abstract}

\section{Introduction}

Let ${\cal D}$ be the space of infinitely smooth functions $t:\bkR^n \to \bkC$ of compact support. The space of Schwartz distributions ${\cal D}'$ is the topological dual of ${\cal D}$, i.e. the space of linear and continuous functionals $F:{\cal D} \to \bkC;\, t \to <F,t>$. The space ${\cal D}'$ provides, in many aspects, a suitable generalization of the set of continuous functions $C^0$.
One of its most interesting properties is that (in a suitable sense) $C^0 \subset {\cal D}'$ while differentiation is an internal operation in ${\cal D}'$ \cite{Zemanian}.
On the other hand, the major limitation of ${\cal D}'$ is that it only displays the structure of a vector space and not of an algebra.
This is has been known since 1954 when L.Schwartz proved that \cite{Schwartz} there is {\it no} associative commutative algebra $({\cal A},+,\circ)$ satisfying the three following properties: (1) The space of distributions ${\cal D}'$ over $\bkR^n$ is linearly embedded into ${\cal A}$ and $f(x)\equiv 1$ is the unity in ${\cal A}$. (2)  The restriction of the product $\circ$ to the set of continuous functions $C^0$ reproduces the pointwise product of functions.
(3) There exist linear derivative operators $\partial_{x_i}:{\cal A}\to {\cal A}$, $i=1,..,n$ satisfying the Leibnitz rule and such that their restriction to ${\cal D}'$ coincide with the usual distributional derivatives.

The best alternative seems to be the Colombeau product \cite{Colombeau} which is defined in a superset ${\cal G} \supset {\cal D}'$ satisfying the properties (1) and (3) with (2) holding only for $C^{\infty}$. Another interesting possibility stems from a solution of the following problem:\\
\\
{\bf Problem:} {\it Determine associative algebras $({\cal A},+,\star)$ in the situation $C^{\infty} \subset {\cal A} \subset {\cal D}'$
and such that: (i) $f\star g= f g$, for all $f,g \in {\cal A}\cap C^0$. (ii) The derivative operators in ${\cal A}$ are of the form $\partial_{x_i} {\cal A} \to {\cal A}$, $i=1,..,n$ and coincide with the restriction to ${\cal A}$ of the usual distributional derivatives in ${\cal D}'$ and, moreover, satisfy the Leibnitz rule.}\\

In this paper we will consider the one dimensional version of this problem, i.e. the case where ${\cal A} \subset {\cal D}'(\bkR^n)$ and $n=1$. For this case we will construct a solution $({\cal A},+,\star)$, explicitly. We will then use the new ${\cal A}$-setting to derive a global formulation for $1$-dimensional quantum confined systems. The details and proofs of the main results will be presented elsewhere \cite{Dias1}. Unfortunatly, it is still not completly clear how these results can be generalized for higher dimensions (i.e. for $n >1$). 

\section{A multiplicative product of distributions}

Let then $C_p^{\infty}$ be the space of piecewise smooth functions on $\bkR$, i.e. the space of functions $f:\bkR \to \bkC$ which are infinitely smooth on $\bkR \backslash V_f$ for some finite set $V_f$ and such that $\lim_{x \to x_0^{\pm}} f^{(n)}(x)$ exists and is finite for all $x_0\in V_f$ and all $n\in \bkN_0$. We define ${\cal A}$ as the set of functions in $C_p^{\infty}$ (regarded as distributions) together with all their distributional derivatives. One can easily prove that $F \in {\cal A}$ iff there is a finite set $V_F=$sing supp $F$, a function $f \in C_p^{\infty}$ and a set of distributions $\Delta_w^{(F)}$ with support on $w \in V_F$ such that $F$ can be written as $F=f + \sum_{w \in V_F} \Delta_w^{(F)}$. Notice that each $\Delta_w^{(F)}$ is supported at $\{w\}$ and so is a finite linear combination of Dirac deltas and its derivatives.

The aim now is to introduce a multiplicative product in the space ${\cal A}$. We start by considering the restriction to ${\cal A}$ of a simple product of distributions with non-intersecting singular supports that was introduced by L. H\"ormander in \cite{Hormander}.\\
\\
{\bf Definition 1: H\"ormander product in ${\cal A}$}\\
{\it Let $F,G \in {\cal A}$ be such that $V_F \cap V_G= {\emptyset } $. Let $\{\Omega_w \subset \bkR, w\in V_F \cup V_G\}$ be a finite covering of $\bkR$ by open sets satisfying $w'\notin \Omega_w$, $\forall w\not=w' \in V_F \cup V_G$. Let us also introduce the compact notation $F_w=F_{\Omega_w}$, $G_w=G_{\Omega_w}$ to designate the restrictions of $F$ and $G$ to the set $\Omega_w$.
Then $F\cdot G$ is defined by its restrictions to $\Omega_w$:
\begin{equation}
F\cdot G: \quad (F\cdot G)_{\Omega_w}=F_wG_w \, ,\quad w \in V_F \cup V_G
\end{equation}
where $F_wG_w$ denotes the usual product of a distribution by an infinitely smooth function.}\\

One can prove that the definition of $F \cdot G$ is independent of the particular covering $\{\Omega_w \}$ and yields a well defined Schwartz distribution, uniquely defined by eq.(2.1), \cite{Dias1,Hormander}.

To proceed we consider the general unrestricted case where $F,G \in {\cal A}$:\\
\\
{\bf Definition 2: The product $\star$} \\
{\it Let $F,G \in {\cal A}$. The multiplicative product $\star$ is defined by:
\begin{equation}
F\star G= \lim_{\epsilon \to 0^+}F \cdot G^{\epsilon}
\end{equation}
where $G^{\epsilon}(x)=G(x+\epsilon)$ is the translation of $G$ by $\epsilon $.}\\

This product is well defined for all $F,G \in {\cal A}$. The functional $F\star G$ acts as: $
<F\star G,t>=\lim_{\epsilon \to 0^+} <F\cdot G^{\epsilon},t>$
and its domain is $\D(\bkR)$, i.e. the limit exists for all $t\in \D(\bkR)$. Hence $F\star G$ is a Schwartz distribution.
Its explicit form is given by:
\begin{equation}
F\star G= fg + \sum_{k=1}^N \left[ g_k
\Delta^{(F)}_{x_k}  +f_{k-1} \Delta^{(G)}_{x_k} \right]
\end{equation}
where $\{x_1<x_2...<x_N\}=V_F \cup V_G$, $x_0=-\infty$ and $x_{N+1}=+\infty$. The functions $f_k,g_k\in C^{\infty}$, $k=0,..,N$ are such that their restrictions to $]x_k,x_{k+1}[$ coincide with the restrictions of $f,g$, respectively. We also defined $\Delta^{(F)}_{x_k}=0$ if $x_k \in V_G \backslash V_F$ and $\Delta^{(G)}_{x_k}=0$ if $x_k \in V_F \backslash V_G$.

Let us point out that the $\star$-product is (i) distributive, (ii) associative but (iii) non-commutative,
(iv) it satisfies the Leibnitz rule, (v) it reproduces the pointwise product of functions and (vi) it reproduces the H\"ormander product for distributions with non-intersecting singular supports. Moreover, it can be consistently extended (through definition 2) to the sets ${\cal A}^{(n)}$ of $C_p^{(2n+1)}$-functions and their distributional derivatives up to order $n+1$ for $n\in \bkN_0$. Notice that ${\cal A}^{(\infty)}={\cal A}$ and that eq.(2.3) is also valid for $F,G \in {\cal A}^{(n)}$.

\section{A global formulation of quantum confined systems}

The algebra $({\cal A},+,\star)$ yields a sufficiently general setting to address some interesting problems. One such application is the global formulation of quantum confined systems. Let us consider a $1$-dimensional dynamical system confined to the positive half-line and described by the Hamiltonian $H=\frac{p^2}{2m} + V(x)$. To keep the discussion simple we make $m=1/2$, $\hbar=1$ and $V(x)=0$. In the usual quantum formulation of the system one defines the Hilbert space to be ${\cal H}=L^2(\bkR^+_0)$ and looks for self adjoint (s.a.) realizations of the Hamiltonian operator. This is easily achieved for $H=-\partial_x^2$ by imposing Dirichlet boundary conditions upon its domain: ${\cal D}(H)= \{\psi \in L^2(\bkR^+_0) : \psi \in AC^1(\bkR^+_0) ; \, \psi(0)= 0 \}$ where $AC^n$ is the set of functions with absolutely continuous $n$-order derivative. The well-known problems \cite{Piotr,Dias2} of this (kinematical) approach to confinement are that i) the usual momentum operator is not s.a. in $L^2(\bkR^+_0)$; ii) there are many different s.a. realizations of the Hamiltonian operator and iii) there is no straightforward translation of this approach to some other formulations of quantum mechanics (e.g. the deformation and the de Broglie-Bohm formulations) \cite{Dias3}.

The problem we want to address is that of formulating the quantum confined system in the global Hilbert space $L^2(\bkR)$. More precisely, we want to derive a s.a., globally defined Hamiltonian operator that dynamical confines the system to $\bkR^+_0$.
Let us start by considering the unconfined eigenvalue equation:
$-\partial_x^2 \psi_U = E \psi_U$.
The confinement to the positive semi-axis implies that the physical eigenfunctions are $\psi_C(x) =\theta(x) \psi_U(x)$ where $\theta$ is the Heaviside step function. By substituting $\psi_C$ into the free eigenvalue equation we get:
\begin{equation}
-\partial_x^2 \psi_C(x) = -\delta'(x) \psi_U(x) -2 \delta(x) \psi_U'(x) + \theta(x) E \psi_U(x)
\end{equation}
and by using the product $\star$ we may rewrite this equation exclusively in terms of $\psi_C$:
\begin{equation}
-\partial_x^2 \psi_C(x) + \delta'(x)\star  \psi_C(x) +2 \delta(x)\star \psi_C'(x) = E \psi_C(x)
\end{equation}
Let us define the operators ($n \in {\bkN}_0$):
\begin{equation}
\hat{\delta}_+^{(n)} : {\cal A}^{(n)}(\bkR) \longrightarrow {\cal A}(\bkR) ; \, \psi \longrightarrow  \delta^{(n)} \star \psi \qquad \mbox{and} \qquad \hat{\delta}_-^{(n)} : {\cal A}^{(n)}(\bkR) \longrightarrow {\cal A}(\bkR) ; \, \psi \longrightarrow  \psi \star \delta^{(n)}
\end{equation}
in terms of which we may rewrite eq.(3.2) as:
\begin{equation}
H_C \psi = \left[-\partial_x^2 + \hat{\delta}'_+(x) +2 \hat{\delta}_+(x) \partial_x \right] \psi_C(x) = E \psi_C(x)
\end{equation}
where we defined the new {\it confined} Hamiltonian $H_C$. One can prove that \cite{Dias2}:
(i) The solutions of $H_C \psi =E \psi$ are exactly $\psi(x) = \theta (x) \psi_U(x) $ where $\psi_U$ is a solution of the free equation.
(ii) The maximal domain of the operator $H_C$ is ${\cal D}_{max}(H_C)= \{\psi \in L^2(\bkR) : \psi(x) = \theta (-x) \psi_-(x) + \theta(x) \psi_+(x); \, \psi_-,\psi_+ \in AC^1(\bkR) ; \, \psi_-(0)= \psi_-'(0)=0 \}$. (iii) $ H_C$ is not symmetric in ${\cal D}_{max}(H_C)$.

To refine our approach in order to derive a s.a. realization of the confined Hamiltonian we may notice that $H_C$ can be written in the form:
\begin{equation}
H_C: {\cal D}_{max}(H_C) \longrightarrow L^2(\bkR); \, \psi \to -\theta(-x) \psi_-''(x) -\theta(x) \psi_+''(x)
\end{equation}
and that this operator displays the symmetric restriction: $H_S \psi=H_C \psi$ for $\psi \in {\cal D}(H_S)=\{\psi \in {\cal D}_{max}(H_C): \psi_+(0)=\psi_+'(0)=0\}$.
This symmetric operator displays a 4-real parameter family of s.a. extensions. The one associated with Dirichlet boundary conditions can be written as:
\begin{equation}
H_D=-\partial_x^2 + \hat{\delta}'_-(x) + \hat{\delta}_-(x) - \hat{\delta}'_+(x) + \hat{\delta}_+(x)
\end{equation}
This operator is s.a. in its maximal domain and furthermore it commutes, in ${\cal D}_{max}(H_D)$, with the projector operators $P_{\pm} \psi=\theta(\pm x) \psi$, i.e. $[H_D,P_{\pm}]\psi=0, \, \forall \psi \in {\cal D}_{max}(H_D)$.
Hence, the solutions of the energy eigenvalue equation can be diagonalized in the representation of $P_{\pm}$. The simultaneous eigenfunctions of $H_D$ and $P_{\pm}$ are $\psi^{(\pm)}_E(x)=\theta(\pm x) \sin kx $ for $k=\sqrt{E}$ and so, in this representation, all the spectral projectors of the Hamiltonian $H_D$ commute with $P_{\pm}$. We conclude that a wave function, originally confined to the positive (or negative) half-line, will stay so forever. Hence, for this simple system, we have achieved a dynamical formulation of quantum confinement. Further developments and detail proofs on this subject have been presented in \cite{Dias2}. Applications to the deformation quantization of confined systems have been discussed in \cite{Dias3}.

\section*{Acknowledgement}
We would like to thank Carlos Sarrico for useful comments and insights into the subject. Research was supported by the grants POCTI/0208/2003 and PTDC/MAT/69635/2006 of the Portuguese Science Foundation.

\end{document}